\documentclass[12pt,bezier]{article}
\usepackage{amsmath}
\usepackage{amsfonts,amsthm,amssymb}
\usepackage{amsfonts}
\usepackage{graphics}
\newtheorem{thm}{Theorem}[section]

\newtheorem{lem}{Lemma}[section]

\newtheorem{conjecture}{Conjecture}
\newtheorem{corollary}{Corollary}[section]

\title{A kind of conditional connectivity of transposition networks generated by $k$-trees
\footnote{The research is supported by NSFC (No.11671296), The Project-sponsored by SRF for ROCS, SEM, and Fund Program for the Scientific Activities of Selected Returned Overseas Professionals in Shanxi Province.}}
\author{ Weihua Yang\footnote{Corresponding author. E-mail: ywh222@163.com.}
\\
\small Department of Mathematics, Taiyuan University of Technology,
Shanxi Taiyuan-030024, China}
\date{}
\usepackage{color}
\begin{document}
\maketitle {\bf Abstract:} {\small  For a graph $G = (V, E)$, a subset $F\subset V(G)$ is called an $R_k$-vertex-cut of $G$ if $G -F$ is disconnected
and each vertex $u \in V(G)- F$ has at least $k$ neighbors in $G -F$. The $R_k$-vertex-connectivity
of $G$, denoted by $\kappa^k(G)$, is the cardinality of the minimum $R_k$-vertex-cut of $G$, which is
a refined measure for the fault tolerance of network $G$. In this paper, we study $\kappa^2$ for Cayley
graphs generated by $k$-trees. Let $Sym(n)$ be the symmetric group on $\{1, 2, \cdots ,n\}$
and $\mathcal{T}$ be a set of transpositions of $Sym(n)$. Let $G(\mathcal{T})$ be the graph on $n$ vertices $\{1, 2, . . . ,n\}$ such that there is an edge $ij$ in $G(\mathcal{T})$ if and only if
the transposition $ij\in \mathcal{T}$. The graph $G(\mathcal{T})$ is called the transposition generating graph of $\mathcal{T}$. We denote by $Cay(Sym(n),\mathcal{T})$ the Cayley graph generated by $G(\mathcal{T})$. The Cayley graph $Cay(Sym(n),\mathcal{T})$ is denoted by $T_kG_n$ if $G(\mathcal{T})$ is a $k$-tree.  We determine $\kappa^2(T_kG_n)$ in this work. The trees are $1$-trees, and the complete graph on $n$ vertices is a $n-1$-tree. Thus, in this sense, this work is a generalization of the such results on Cayley graphs generated by transposition generating trees\cite{yang2} and the complete-transposition graphs\cite{wang}.  }

 {\flushleft\bf Keywords}: Star graphs; Bubble sort graphs; Transposition networks;
  Fault-tolerance; Transposition generating trees

 \section{Introduction}

The interconnection network of a communication or distributed
computer system is usually modeled by a (directed) graph $G$ in which
the vertices represent the switching elements or processors and
communication links are represented by (directed) edges. The traditional  connectivity $\kappa$ and the edge connectivity $\lambda$ of a network
are two classic parameters measuring fault tolerance. The higher these parameters are, the more
reliable the network is [16].  However, it always
underestimates the resilience of large networks. There is a
discrepancy because the occurrence of events which would disrupt a
large network after a few processor or link failures is highly
unlikely. Thus the disruption envisaged occurs in a worst case
scenario (see \cite{esfa,latifi}  for a detailed explanation for the
shortcoming of using $\kappa(G)$ to measure the network
reliability). To overcome the shortcoming, Esfahanian~\cite{esfa}
proposed the concept of restricted connectivity which is a special
case of conditional connectivity proposed by Harary \cite{harary}.
This concept was generalized by Latifi et al. \cite{latifi} to
$R_k$-vertex-connectivity as a measure of conditional fault
tolerance of networks.

A $R_k${\em -vertex-set} of a graph $G$ is a vertex subset
$F\subseteq V(G)$  such that every vertex $u \in V(G) - F$ has at
least $k$ neighbors in $G - F$.  A $R_g${\em-vertex-cut} of a connected graph $G$ is a
$R_k$-vertex-set $F$ such that $G - F$ is disconnected. The
cardinality of a minimum $R_k$-cut of $G$ is the
{\em$R_k$-vertex-connectivity} of $G$, denoted by $\kappa^k(G)$.

Let $\Gamma$ be a group and $S$ be a subset of $\Gamma\setminus \{1_{\Gamma}\}$, where $1_{\Gamma}$ is the identity of $\Gamma$. Cayley digraph $Cay(\Gamma, S)$ is the digraph with
vertex set $\Gamma$ and arc set $\{(g, g\cdot s): g\in\Gamma, s \in S\}$. We say that arc $(g, g\cdot s)$ has label $s$. In particular, if $S^{-1} = S$, then $Cay(\Gamma, S)$ is an
undirected graph, called Cayley graph.

Cayley (di) graphs have a lot of properties which are desirable in an interconnection network \cite{heydemann,jwo}: vertex symmetry
makes it possible to use the same routing protocols and communication schemes at all nodes; hierarchical structure facilitates
recursive constructions; high fault tolerance implies robustness, among others. The transposition networks, Cayley graphs generated by transpositions,  is the most popular family such Cayley graphs. In particular, star graphs~\cite{akers} and  bubble sort graphs~\cite{chou} are the most famous ones, see the survey \cite{heydemann,jwo,latifi2} for other kinds of Cayley graphs.

Since $\kappa(G)\leq \lambda(G)\leq \delta(G)$, a graph $G$ with $\kappa(G)=\delta(G)$ is called maximally vertex connected and a
graph G with $\lambda(G)=\delta(G)$ is called maximally edge connected. It is well-known~\cite{tindell} that any vertex transitive graph is maximally
edge connected (Cayley (di) graphs are always vertex transitive), and many Cayley graphs are maximally vertex connected
( such as hypercubes, star graphs, bubble sort graphs, etc.). Works on $\kappa^1$ of  Cayley graphs generated by transposition trees can be found in \cite{eddie1,eddie2,huyang}. In \cite{wan}, Wan and Zhang determined $\kappa^2$ for star graphs, which was generalized by Cheng et al. \cite{eddie3}  to Cayley graphs generated by transposition trees (independently by Yang~\cite{yang2}). Zhang et al. in \cite{zhang1} determined $\kappa^2$ for alternating group graphs, and Cheng et al. \cite{eddie4} generalized the results in \cite{zhang1} to the Cayley graphs generated by 2-trees. Recently, Yu et al. \cite{zhang2}  determined $\kappa^2$ for Cayley graphs generated
by unicyclic graphs. Wang et al. \cite{wang} determined $\kappa^2$ for Cayley graphs generated
by complete transposition graphs. For the related research on the $R_k$-vertex connectivity for the Cayley graphs, we refer to \cite{oh,wu,yang,yang1,yang3}.

In this work, we consider Cayley graph $Cay(Sym(n),\mathcal{T})$, where $Sym(n)$ is the symmetric group on $\{1, 2,. . . ,n\}$ and $\mathcal{T}$ is a set
of transpositions of $Sym(n)$. Let $G(\mathcal{T})$ be the graph on $n$ vertices $\{1, 2, . . . ,n\}$ such that there is an edge $ij$ in $G(\mathcal{T})$ if and only if
the transposition $ij\in \mathcal{T}$. The graph $G(\mathcal{T})$ is called the transposition generating graph of $\mathcal{T}$. For convenience, we call $Cay(Sym(n),\mathcal{T})$ the Cayley graph generated by $G(\mathcal{T})$.  In particular, if $G(\mathcal{T})$  is a path,
$Cay(Sym(n),\mathcal{T})$ is the  bubble-sort graph, denoted by $BS_n$. If $G(\mathcal{T})$  is a star, $Cay(Sym(n),\mathcal{T})$ is the  well-known star graph, denoted
by $S_n$.  If $G(\mathcal{T})$  is a tree, the $Cay(Sym(n),\mathcal{T})$ includes the $BS_n$ and $S_n$ as its subclasses.

 A $k$-tree $T_{k,n}$ with $n$ vertices is defined recursively as follows: A set of $k + 1$ mutually adjacent vertices constitutes a $k$-tree $T_{k,k+1}$ and a $k$-tree $T_{k,n+1}$ is any graph obtained by joining a new vertex to $k$ mutually adjacent vertices of a $k$-tree $T_{k,n}$. One can see that if $k=1$, then $T_{k,n}$ is a tree. Thus, the $k$-trees are a generalization of trees, such a generalization of trees has been studied extensively, first in \cite{beineke1,beineke2}. If $G(\mathcal{T})$ is a $k$-tree,  $Cay(Sym(n),\mathcal{T})$ is denoted by $T_kG_n$. In particular, if  $k=n-1$, then $T_{n-1,n}$ is the complete graph. If $G(\mathcal{T})$ is complete, $Cay(Sym(n),\mathcal{T})$ is called complete-transposition graphs, denoted by $CT_n$. Wang et al. determined $\kappa^2(CT_n)$ in \cite{wang}. In this paper, we  completely determined $\kappa^2(T_kG_n)$, which generalizes the results on Cayley graphs generated by transposition trees and by complete transposition graphs, as  transposition trees are $T_{1,n}$ and the complete transposition graphs are $T_{n-1,n}$.

 \begin{center}
\scalebox{0.6}{\includegraphics{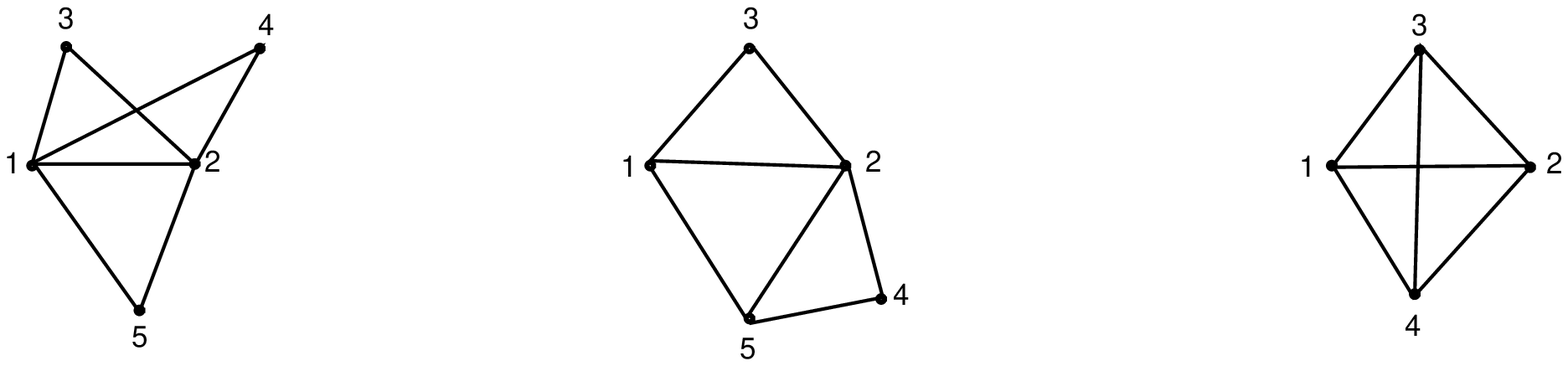}}\\
Figure~1. 2-trees on 5 vertices and 3-tree on 4 vertices
\end{center}

 \section{Preliminaries}

For a graph $G$, a subgraph $G_1$ of $G$, and a vertex $u \in V(G)$, the neighbor set of $u$ in $G_1$ is $N_{G_1} (u) =\{v $ is adjacent with
$u$ in $G\}$. In particular, if $G_1 = G$, then $N_G(u)$ is the neighbor set of $u$ in $G$ and $d_G(u) = |N_G(u)|$ is the degree of vertex $u$ in $G$. The
minimum degree of $G$ is denoted as $\delta (G)$. For a subset $U\in V(G), N_{G}(U)=\cup_{u\in U} N_G(u)-U$, and $G[U]$ is the subgraph of $G$ induced
by $U$. Sometimes, we use a graph itself to represent its vertex set. For example, $N_G(G_1)$ is used to denote $N_G(V(G_1))$
where $G_1$ is a subgraph of $G$, and $N_F(U)$ is used to denote $N_{G[F]}(U)$ where $F, U$ are vertex sets. When graph $G$ is obvious in
the context, we omit the subscript $G$ and use $N(U)$ to denote $N_G(U)$. A cycle with length $k$ is called a $k$-cycle. The length of
the shortest cycle of $G$ is called the girth of $G$, denoted by $g(G)$.

 In this section, we first introduce some notations and former results which will be used in our proofs. Then we shall derive
some structural properties for  $T_kG_n$. By the definition, a $k$-tree $T_{k,n}$ has $kn-\frac{k(k+1)}2$ edges and $T_{k,n}$ contains a triangle if $k\geq 2$. Thus $T_kG_n$ is a  $kn-\frac{k(k+1)}2$ regular graph on $n!$ vertices. For convenience, we assume $n$ is the last vertex added on $G(\mathcal{T})$, that is, the degree of $n$ in  $G(\mathcal{T})$ is $k$. And we call a vertex of degree $k$ in $T_{k,n}$ is a leaf. Moreover, we assume $k\geq2$, as the case for $k=1$ has been considered in \cite{eddie3,yang2}.

 \begin{center}
\scalebox{0.8}{\includegraphics{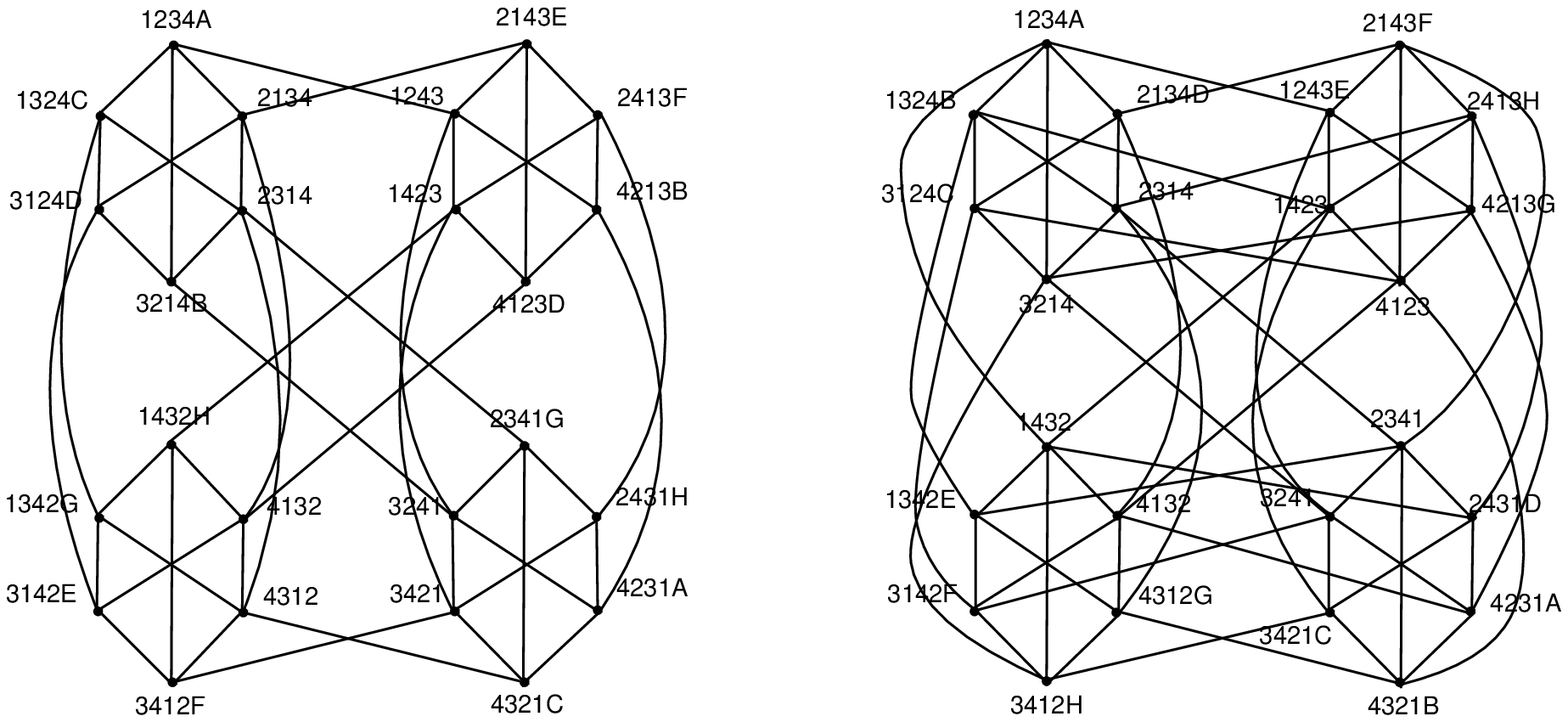}}\\
Figure~2. Cayley graphs $T_2G_4$ and $T_3G_4$.
\end{center}

 One most important property of transposition networks is the hierarchical structure, which is particularly useful
in our inductive arguments.

Assume $k<n-1$ and $n$ is a leaf of $G(\mathcal{T})$. One can see that
$T_kG_n$  can be decomposed into $n$ interconnected copies of
$T_kG_{n-1}$, say $T_kG_{n-1}^1,T_kG_{n-1}^2,
\cdots,T_kG_{n-1}^n$, where $T_kG_{n-1}^i$ is the  subgraph  induced
by vertex set $\{(p_1p_2\cdots p_{n-1}i)~|(p_1p_2\cdots
p_{n-1})$ ranges over all permutations of $\{1,2,\cdots,
n\}\setminus \{i\}$. In particular, if $k=n-1$, then  $T_kG_n$  can be decomposed into $n$ interconnected copies of
$T_{k-1}G_{n-1}$, say $T_{k-1}G_{n-1}^1,T_{k-1}G_{n-1}^2,
\cdots,T_{k-1}G_{n-1}^n$. Without loss of the generality,  we always decompose $T_kG_n$ into $n$ copies by using the leaf $n$. For a vertex $u\in V(T_kG_{n-1}^i)$, we call the neighbors out of $T_kG_{n-1}^i$ outsider neighbors. By the arguments above, one can see the following.

\begin{lem}\label{property}
Let $1\leq i,j\leq n$ be two distinct integers. Then the following holds.
\begin{enumerate}
\item Any vertex $u\in V(T_kG_{n-1}^i)$ has $k$ outsider neighbors in $k$ different copies.

\item The outsider neighbors of the vertices in $T_kG_{n-1}^i$ are all different.

\item There are $k(n-2)!$ independent edges between  $T_kG_{n-1}^i$ and $T_kG_{n-1}^j$. That is, $|N(T_kG_{n-1}^i)\cap V(T_kG_{n-1}^j)|=k(n-2)!$
\end{enumerate}
\end{lem}
{\flushleft\textbf{Proof.}}\quad Suppose $u=p_1p_2\cdots p_{n-1}i\in V(T_kG_{n-1}^i)$ and let $i_1,i_2,\cdots,i_k$ be the $k$ neighbors of $n$ in   $G(\mathcal{T})$. Then $u$ has $k$ outsider neighbors   $u(i_1n)\in V(T_kG_{n-1}^{p_{i_1}})$, $u(i_2n)\in V(T_kG_{n-1}^{p_{i_2}})$, $\cdots$, $u(i_kn)\in V(T_kG_{n-1}^{p_{i_k}})$. Notice that $p_{i_k}\not=p_{i_l}$ for $i_k,i_l\in\{i_1,i_2,\cdots,i_k\}$. Thus, $(1)$  holds. Since $(1)$  holds, $(2)$ is clearly true. For a vertex  $u=p_1p_2\cdots p_{n-1}i\in V(T_kG_{n-1}^i)$, if $j\in\{p_{i_1},\cdots,p_{i_k}\}$, then there is an outsider neighbor of $u$ in $T_kG_{n-1}^j$. Assume $p_{i_1}=j$, then there are $(n-2)!$ such vertices in $T_kG_{n-1}^i$ and then there are $(n-2)!$ such outsider neighbors in  $T_kG_{n-1}^j$. So the subgraph $T_kG_{n-1}^i$ implies $k(n-2)!$ outsider neighbor in  $T_kG_{n-1}^j$. Therefore, $(3)$ holds.\hfill$\Box$

The following result is due to Cheng and Lipt\'ak \cite{eddie1} on the Cayley graphs  generated by transpositions.

\begin{thm}
Let $G$ be the Cayley graph obtained from a generating graph $G(\mathcal{T})$ on $\{1, 2, \cdots , n\}$ with $m$ edges, where $m\geq 7$. Suppose
that $T$ is a set of vertices of $G$ such that $|T|\leq 2m -2$ if $G(\mathcal{T})$ has no triangles and  $|T|\leq 2m -3$ if $G(\mathcal{T})$  has a triangle. Then one of the
following is true:
\begin{enumerate}\label{eddie}
\item $G- T$ is connected.

\item $G -T$ is disconnected with exactly two components, one of which is a singleton.

\item The transposition generating graph $G(\mathcal{T})$ has no triangles, $G -T$ is disconnected with exactly two components, one of which is $K_2$, and
$|T| = 2m -2$.

\item The transposition generating graph $G(\mathcal{T})$ has no triangles, $G -T$ is disconnected with exactly three components, two of which are singletons, and
$|T| = 2m -2$.

\item The transposition generating graph $G(\mathcal{T})$ has a triangle, $G -T$ is disconnected with exactly three components, two of which are singletons, and
$|T| = 2m -3$.
\end{enumerate}
\end{thm}

Notice that a $k$-tree has  $kn-\frac{k(k+1)}2$ edges. By using the lemma above, we have the following without proof.

\begin{corollary}
For any integer $n\geq5,2\leq k\leq n-1$, let $T_kG_n$ be the Cayley graphs generated by a $k$-tree. Then $\kappa^1(T_kG_n)=2(kn-\frac{k(k+1)}2)-2$.

\end{corollary}

The girth of $T_kG_n$ is 4. For a 4-cycle
$(u_1u_2u_3u_4)$ where $u_2 = u_1(ij), u_3 = u_2(kl), u4 = u_3(ij), u_1 = u_4(kl)$, and $i, j, k, l$ are all distinct. Call this form of 4-cycle as Type A 4-cycle.
When $G(\mathcal{T})$ contains triangles, except Type A 4-cycles, there are 4-cycles in $T_kG_k$ having the form $(u_1u_2u_3u_4)$ where $u_2 = u_1(sp), u_3 = u_2(pt), u_4 = u_3(sp),
u_1 = u_4(st)$, and $spt$ is the unique 3-cycle in $G(\mathcal{T})$. Call this form of 4-cycle as Type B 4-cycle.

For any two distinct vertices $u, v$ in $T_kG_n$,  $|N(u) \ N(v)| \leq 3$. Furthermore, if $|N(u) \ N(v)| = 3 $ , then the three common
neighbors of $u$ and $v$ have the form $u_1 = u(st) = v(sp), u_2 = u(sp) = v(pt)$, and $u_3 = u(pt) = v(st)$, where $spt$ is a triangle in
$G(\mathcal{T})$. The Cayley graph $T_kG_n$ has an important structural property bolew.

\begin{lem}[\cite{eddie1}]\label{k24}
The Cayley graph $T_kG_n$ contains no $K_{2,4}$ as subgraphs.
\end{lem}

Recall that we denote the grith of a graph $G$ by $g(G)$. We list the following known results due to Wang and Zhang.

\begin{thm}[\cite{wan}]\label{star}
For any integer $n\geq4$,  $\kappa^2(S_n)=g(n-3)$.
\end{thm}

Cheng and Lipt\'ak~\cite{eddie3}, and Yang and Meng~\cite{yang2} generalized the result above to the Cayley graphs generated by transposition generating trees.

\begin{thm}[\cite{eddie3,yang2}]\label{tree}
Let $G_n$ be the Cayley graph generated by transposition generating tree. Then  $\kappa^2(G_n)=g(n-3)$ for $n\geq4$.
\end{thm}

Recently, Wang et al.~\cite{wang} consider the complete transposition graphs.

\begin{thm}[\cite{wang}]\label{complete}
Let $CT_n$ be the complete transposition graph. Then  $\kappa^2(CT_n)=2n(n-1)-10$ for $n\geq5$.
\end{thm}

For the Cayley graph generated by unicyclic graphs $UC_n$ (an unicyclic graph means the graph contains exactly one cycle), Yu et al. showed the following.

\begin{thm}[\cite{zhang2}]\label{uc}
Let $UC_n$ be the Cayley graph generated by a unicyclic graph $G(\mathcal{T})$ and $n\geq4$. Then  $\kappa^2(UC_n)=4n-10$ if $G(\mathcal{T})$ contains a triangle, and  $\kappa^2(UC_n)=4n-10$ if  $G(\mathcal{T})$ contains no  triangles.
\end{thm}

\section{Determining $\kappa^2(T_kG_n)$}

In this section, we shall show that $\kappa^2(T_kG_n)=4(kn-\frac{k(k+1)}2)-10$ for $n\geq5,k\geq2$.  Wang et al.~\cite{wang} showed that $\kappa^2(T_3G_4)=16$, but $16\not=4(kn-\frac{k(k+1)}2)-10$. So $n\geq5$ is necessary.


\begin{thm}\label{main}
For any integers $n\geq 5,k\geq2$, $\kappa^2(T_kG_n)=4(kn-\frac{k(k+1)}2)-10$.
\end{thm}

{\flushleft\textbf{Proof.}}\quad If $k=n-1$, then the claim holds by the main results in \cite{wang}. From now on, we assume $2\leq k\leq n-2$.

{\flushleft\textbf{Claim.}}\quad $\kappa^2(T_kG_n)\leq 4(kn-\frac{k(k+1)}2)-10$.

Assume $ijk$ be a triangle in $\mathcal{G(T)}$. Then $u,v=u(ij),w=v(jk),z=u(ik)=w(ij)$ forms a type B 4-cycle in $T_kG_n$, denoted by $C$. Notice that $x=u(jk)=w(ik),y=v(ik)=z(jk)$ and $T_kG_n$ does not contains the $K_{2,4}$ as subgraph. Then we have $|N(C)|=4(kn-\frac{k(k+1)}2)-10$.  We shall show that $F=N(C)$ is an $R_2$-vertex-cut of $T_kG_n$. So it suffices to show
that every vertex in $T_kG_n-F$ has at least two  neighbors. This is true for vertices in $C$. Let $t$ be a vertex in $T_kG_n-F-V(C)$.
Since $T_kG_n$ is a bipartite graph and thus contains no odd cycle, we must have $N_F(t)\subset N(\{u,w\})$ or $N_F(t)\subset N(\{v,z\})$. Suppose,
without loss of generality, that $N_F(t)\subset N(\{u,w\})$. By Lemma~\ref{k24}, $u$ and $t$ have at most 3 common neighbors,
otherwise, constitute a $K_{2,4}$. Similarly, $w$ and $t$ have at most 3 common neighbors. That is, $t$ is adjacent to at most six vertices
in $F$. So, if $n\geq 6$, or $n=5$ and $k= 3$, then $t$ has at least $kn-\frac{k(k+1)}2-6\geq 3$ neighbors in $T_kG_n-F-V(C)$. That is, $F$ is an $R_2$-vertex-cut of $T_kG_n$ for the cases.
Suppose $n=5,k=2$. The Cayley graph $T_kG_n$ can be decomposed in $5$ copies of $T_2G_4$ (see Figure 2 for the $T_2G_4$). We pack the 4-cycle $C$ formed by ${u=1324,v=2314,w=2134,z=3124}$. One can see that $F$ induces an  $R_2$-vertex-cut of the $T_2G_4$ including the 4-cycle  (in $T_kG_n$, say $T_2G_4^1$). Then, if $t$ is in $T_2G_4^1$, then $t$ has at least two neighbors in $T_kG_n-F-V(C)$. By  Lemma~\ref{property} (1), if $t\not\in V(T_2G_4^1)$, then $t$ has at most one neighbor in $F\cap V(T_2G_4^1)$. And $\{u,w\}$ ($\{v,z\}$) has four outsider neighbors. So $t$ has at most five neighbors in $N(\{u,w\})$. That is, $t$ has at least $kn-\frac{k(k+1)}2-5=7-5=2$ neighbors in $T_kG_n-F-V(C)$. Thus, $F=N(C)$ is an $R_2$-vertex-cut of $T_kG_n$ with $4(kn-\frac{k(k+1)}2)-10$ vertices. The claim holds.

We next show that $\kappa^2(T_kG_n)\geq 4(kn-\frac{k(k+1)}2)-10$. By contradiction, we suppose that $F$ is an $R_2$-vertex-cut with size no more than $4(kn-\frac{k(k+1)}2)-11$. Let $F_i=V(T_kG_{n-1}^i)\cap F$.

{\flushleft\textbf{Case 1.}}\quad $|F_i|\leq 2(k(n-1)-\frac{k(k+1)}2)-4$ for all $i$.

Notice that $T_kG_n$ is maximally connected. Then by Lemma~\ref{eddie},  $T_kG_{n-1}^i-F_i$ either is connected or has at most one singleton. Let $J=\{i:T_kG_{n-1}^i-F_i$ has a singleton$\}$. A simple count shows that $|J|\leq 4$.

{\flushleft\it{Subcase 1.1}}\quad  $T_kG_{n-1}^i-F_i$ is connected for all $i$.

By Lemma~\ref{property}, there are $k(n-2)!$ independent edges between $T_kG_{n-1}^i$ and $T_kG_{n-1}^j$ for any pair $\{i,j\}$. Notice that $T_kG_{n-1}^i-F_i\subset F_j$ if $T_kG_{n-1}^i-F_i$ is disconnected to  $T_kG_{n-1}^j-F_j$. One can see that the inequality $k(n-2)!>2[2(k(n-1)-\frac{k(k+1)}2)-4]$ holds if $n\geq6$.  Thus, $T_kG_{n-1}^i-F_i$ is connected to  $T_kG_{n-1}^j-F_j$ for all $j\not=i$, that is, $T_kG_{n}-F$ is connected if $n\geq 6$, or $n=5$ and $k=3$, a contradiction. When $n=5$ and $k=2$, the inequality becomes an equality. That is, if there exists  one pair  $\{i,j\}$  such that $T_kG_{n-1}^i-F_i$ is disconnected to  $T_kG_{n-1}^j-F_j$, then $|F_i|=|F_j|=2(k(n-1)-\frac{k(k+1)}2)-4$ and $|F|-|F_i|-|F_j|=2$. Therefore, there must exist an integer $l$ such that $F_l=\emptyset$. Since $k(n-2)!>2(k(n-1)-\frac{k(k+1)}2)-4\geq |F_i|$ for $n\geq5$, $T_kG_{n-1}^l$ is connected to  $T_kG_{n-1}^j-F_j$ for each $j$. Thus, $T_kG_{n}-F$ is connected if $n\geq 5$, a contradiction.

{\flushleft\it{Subcase 1.2}}\quad There is exactly one $i$ such that $T_kG_{n-1}^i-F_i$ is disconnected.

Without loss of generality, we assume  $T_kG_{n-1}^1-F_1$ is disconnected, that is,  $T_kG_{n-1}^1-F_1$ contains exactly two components such that one of them is a singleton $v$. Since $T_kG_n$ is maximally connected, we have $|F_1|\geq k(n-1)-\frac{k(k+1)}2$. By an argument similar to that of Subcase 1.1, one can see that  $T_kG_{n-1}^2-F_2,\cdots,T_kG_{n-1}^n-F_n$ are in the same component $C$ of $T_kG_n-F$. Now we show that the large component $H$ of $T_kG_{n-1}^1-F_1$ is also in $C$. Notice that $T_kG_{n-1}^1-F_1-\{v\}$ has  $k[(n-1)!-|F_1|-1]$ outsider neighbors. However, $|F-F_1|<k[(n-1)!-|F_1|-1]$. So $T_kG_{n-1}^1-F_1$ is  in the large component $C$, and $\{v\}$ is either a component or in $C$. This contradicts with $F$ being an $R_2$-vertex-cut.

{\flushleft\it{Subcase 1.3}}\quad There is exactly two integer $i$'s  such that $T_kG_{n-1}^i-F_i$ is disconnected.

Without loss of generality, we assume that $T_kG_{n-1}^1-F_1$ and $T_kG_{n-1}^2-F_2$ are disconnected. Let $v_1$ and $v_2$ be the singletons of $T_kG_{n-1}^1-F_1$ and $T_kG_{n-1}^2-F_2$ respectively. Thus, we have $|F_1|\geq k(n-1)-\frac{k(k+1)}2$ and $|F_2|\geq k(n-1)-\frac{k(k+1)}2$. Similarly, $T_kG_{n-1}^3-F_3,\cdots,T_kG_{n-1}^n-F_n$ are in the same component $C$ of $T_kG_n-F$. Since $T_kG_{n-1}^1-F_1-\{v_1\}$ ($T_kG_{n-1}^2-F_2-\{v_2\}$) has at least $k[(n-1)!-|F_1|-1]-k(n-2)!$ outsider neighbors in $T_kG_n-V(T_kG_{n-1}^2)$ and  $|F-F_1-F_2|<k[(n-1)!-|F_1|-1]-k(n-2)!$, we have that the larger components in $T_kG_{n-1}^1-F_1$  ($T_kG_{n-1}^2-F_2$) is also in the same component $C$. The subgraph induced by $\{v_1,v_2\}$ has degree at most one. So this contradicts with $F$ being an $R_2$-vertex-cut again.

{\flushleft\it{Subcase 1.4}}\quad There is exactly three integer $i$'s  such that $T_kG_{n-1}^i-F_i$ is disconnected.

Without loss of generality, we assume that $T_kG_{n-1}^1-F_1$, $T_kG_{n-1}^2-F_2$ and  $T_kG_{n-1}^3-F_3$ are disconnected such that $v_1,v_2$ and $v_3$ are the corresponding singletons, respectively. Similarly,  $T_kG_{n-1}^4-F_4,\cdots,T_kG_{n-1}^n-F_n$ are in the same component $C$ of $T_kG_n-F$. We show that $T_kG_{n-1}^i-F_1-\{v_i\}, i=1,2,3$  is also in $C$. Notice that $T_kG_{n-1}^i-F_1-\{v_i\},( i=1,2,3)$ has at least $k[(n-1)!-|F_1|-1]-2k(n-2)!$ outsider neighbors in $T_kG_n-V(T_kG_{n-1}^2\cup T_kG_{n-1}^3)$ and $k[(n-1)!-|F_1|-1]-2k(n-2)!>|F-F_1-F_2-F_3|$. So $T_kG_{n-1}^i-F_1-\{v_i\}, i=1,2,3$  is also in $C$. Clearly, the subgraph induced by $\{v_1,v_2,v_3\}$ contains a vertex of degree at most one, a contradiction.

{\flushleft\it{Subcase 1.5}}\quad There is exactly four integer $i$'s  such that $T_kG_{n-1}^i-F_i$ is disconnected.

Without loss of generality, we assume that $T_kG_{n-1}^1-F_1$, $T_kG_{n-1}^2-F_2$, $T_kG_{n-1}^3-F_3$ and $T_kG_{n-1}^4-F_4$  are disconnected such that $v_1,v_2$, $v_3$ and  $v_4$ are the corresponding singletons, respectively. Similarly, $T_kG_{n-1}^i-F_i-\{v_i\}$ ($i=1,2,3,4$) and $T_kG_{n-1}^i-F_i-\{v_i\}$ ($i\geq 5$) are in the same component of $T_kG_n-F$. That is, if $F$ is an  $R_2$-vertex-cut, then $\{v_1,v_2,v_3,v_4\}$ forms a 4-cycle.

Without loss of generality, we suppose that $v_1v_2v_3v_4v_1$ is a 4-cycle and $v_1=p_1p_2\cdots$ $2\cdots 3\cdots$ $4\cdots p_{n-1}1$. Then we have $v_2=p_1p_2\cdots 1$ $\cdots 3$ $\cdots 4$ $\cdots p_{n-1}2$ and $v_4=p_1p_2\cdots2$ $\cdots 3$ $\cdots 1$ $\cdots p_{n-1}4$. On the one hand, $v_3=p_1p_2\cdots 1$ $\cdots 2$ $\cdots 4$ $\cdots p_{n-1}3$ since $v_2$ is a neighbor of $v_3$ by using some edge due to the leaf $n$. On the other hand, $v_3=p_1p_2\cdots 2$ $\cdots 4$ $\cdots 1$ $\cdots p_{n-1}3$ since $v_4$ is a neighbor of $v_3$. This is impossible.

{\flushleft\textbf{Case 2.}}\quad $|F_i|\geq 2(k(n-1)-\frac{k(k+1)}2)-3$ for some $i$.

Let $I=\{i:|F_i|\geq 2(k(n-1)-\frac{k(k+1)}2)-3\}$. Since $|F|\leq 4(kn-\frac{k(k+1)}2)-11$, we have $|I|\leq 2$. For $j\not\in I$, $T_kG_{n-1}^j-F_j$ contains at most one singleton $u_j$ (if it exists).

{\flushleft\it{Subcase 2.1}}\quad $|I|=1$.

Without loss of generality, we assume $I=\{1\}$. We first assume that  $C=\cup_{j\not=1}T_kG_{n-1}^j-F_j$ is connected. Since $F$ is an $R_2$-vertex-cut, there is a component $H$ of $T_kG_n-F$ in $T_kG_{n-1}^1-F_1$. Clearly, each vertex of $H$ has at least two neighbors in $T_kG_{n-1}^1-F_1$. Let $P=v_1v_2v_3v_4$ be a path of $H$. Since $H$ is a component of $T_kG_n-F$, each vertex $x\in N_{T_kG_{n-1}^1}(P)$
is either in $F$ or its outsider neighbors in $F$. Thus $|F|\geq |N(P)|\geq 4(kn-\frac{k(k+1)}2)-10$, a contradiction.

Assume $C$ is not connected. Thus there is some $T_kG_{n-1}^j-F_j$ containing two component one which is a singleton. Similarly, the larger component of $T_kG_{n-1}^j-F_j$ for $j=2,3,\cdots$ are in the same component of $T_kG_n-F$, denoted by $C'$. Let $U$ be set of singletons of $T_kG_{n-1}^j-F_j$ which are disconnected to $C'$. Clearly, $U\not=\emptyset$ since $C$ is not connected. Notice that $F$ is an $R_2$-vertex-cut and $u_j$ has at most one outsider neighbor in $T_kG_{n-1}^1$. So $|U|\geq2$. A simple count shows that $|U|\leq3$ since $|F|\leq 4(kn-\frac{k(k+1)}2)-10$. Assume $|U|=3$ and $U=\{u_2,u_3,u_4\}$. Each neighbor of $u_j,j=2,3,4$ is either in $F$ or in $T_kG_{n-1}^1$. Then we have $|F|\geq |F_1|+3(k(n-1)-\frac{k(k+1)}2)+(k-3)+2(k-2)\geq 3(k(n-1)-\frac{k(k+1)}2)-3+2(k(n-1)-\frac{k(k+1)}2)+(k-3)+2(k-2)>4(kn-\frac{k(k+1)}2)-11\geq |F|$, a contradiction. So $|U|=2$.

Assume, without loss of generality, $U=\{u_2,u_3\}$. Since $F$ is an $R_2$-vertex-cut and $u_j$ has at most one neighbor in $T_kG_{n-1}^1$ for $j=2,3$, $u_2$ is a outsider neighbor of $u_3$. Similarly, one can see that the components in $T_kG_{n-1}^1-F_1$, which are disconnected to the edge $u_2u_3$, are connected to $C'$. Let $u_2'$ and $u_3'$ be the outsider neighbors of $u_2$ and $u_3$ respectively. Then $u_2'$ (and $u_3'$) has $k-1$ outsider neighbors in $\cup_{i=4}^nT_kG_{n-1}^i$ by Lemma~\ref{property}. We have known that the edges with label $(in)$ due to the leaf $n$ can not induce a 4-cycle. So the outsider neighbor of $u_2'$ ($u_3'$ ) is not the neighbor of $u_2$ and $u_3$. This implies $u_2'$ and $u_3'$ have $2(k-1)$ outsider neighbors in $F-F_1-F_2-F_3-N_{\cup_{i=4}}^n(u_2)-N_{\cup_{i=4}}^n(u_3)$. Notice that $|F|-|F_1|-|N_{T_kG_{n-1}^2}(u_2)|-|N_{T_kG_{n-1}^3}(u_3)|-|N_{\cup_{i=4}}^n(u_2)|-
|N_{\cup_{i=4}}^n(u_3)|\geq 4(kn-\frac{k(k+1)}2)-11-[2(k(n-1)-\frac{k(k+1)}2)-3]-2[k(n-1)-\frac{k(k+1)}2]-2(k-2)<2(k-1)$. This is impossible.

{\flushleft\it{Subcase 2.2}}\quad $|I|=2$.

Without loss of generality, we assume $I=\{1,2\}$.  There is at most one integer $j$ such that $T_kG_{n-1}^j-F_j$ is disconnected and one of its two component is the singleton $u_j$. Let $C=\cup_{j\not=1,2}T_kG_{n-1}^j-F_j$.

First, we assume $C$ is connected. Then there is a component $H$ of $T_kG_n-F$ in $(T_kG_{n-1}^1-F_1)\cup (T_kG_{n-1}^2-F_2)$ since $F$ is an $R_2$-vertex-cut. Clearly, the degree of each vertex in $H$ is at least 2. Notice that $T_kG_n$ is bipartite and contains no the $K_{2,4}$ as subgraph. A simple count shows that $|N(H)|>|F|$ if $4\leq |V(H)|\leq 6$. This is impossible. Thus we assume $|V(H)|\geq 7$. Then we may take four vertices from $H\cap T_kG_{n-1}^1$ (or $H\cap T_kG_{n-1}^2$), say $v_1,v_2,v_3,v_4$. The four vertices have $4k$ outsider neighbors and at most four of them are in $T_kG_{n-1}^2$. So there are at least $4k-4$ outsider neighbors of $v_1,v_2,v_3,v_4$ in $\cup_{j=3}^nT_kG_{n-1}^j$.  However, $|F|-|F_1|-|F_2|\leq 4k-5$. This implies that $H$ is connected to $C$, a contradiction.

We next assume $C$ is disconnected. Then $C$ consists of two components such that one of them is a singleton. Assume that $u$  is the singleton of $C$ and let $C'=C-\{u\}$. Without loss of generality, we assume $u\in T_kG_{n-1}^3$. Similarly, there is a component $H$ of $T_kG_n-F$ in $(T_kG_{n-1}^1-F_1)\cup (T_kG_{n-1}^2-F_2)\cup \{u\}$. If $u$ is not in $H$, then we are done (see the argument above). Thus, $u\in V(H)$. Note $u$ has at least $k-2$ outsider neighbors in $\cup_{j=4}^n(T_kG_{n-1}^j)$. Similarly, one may take $v_1,v_2,v_3$ from $H\cap T_kG_{n-1}^1$ (or $H\cap T_kG_{n-1}^2$), and $v_1,v_2,v_3$ have at least $3(k-2)$ outsider neighbors in $\cup_{j=4}^n(T_kG_{n-1}^j)$. Thus, $4(k-2)\leq |F-F_1-F_2-F_3|\leq 4k-5-[k(n-1)-\frac{k(k+1)}2]$. This implies $3-[k(n-1)-\frac{k(k+1)}2]>0$, a contradiction. \hfill$\Box$

Combining  Theorem~\ref{star}, Theorem~\ref{tree}, and Theorem~\ref{complete}, we have the following corollary.

\begin{corollary}
Let $T_kG_n$ be a Cayley graph generated by a $k$-tree $G(\mathcal{T})$ and $n\geq 5$.

 $(1)$ If $k\geq2$, then $\kappa^2(T_kG_n)=4(kn-\frac{k(k+1)}2)-10$;

$(2)$  If $k=1$ and $G(\mathcal{T})$ is not a star, then $\kappa^2(T_kG_n)=4(kn-\frac{k(k+1)}2)-8$;

$(3)$ If $G(\mathcal{T})$ is  a star, then $\kappa^2(T_kG_n)=6(kn-\frac{k(k+1)}2)-12$.
\end{corollary}

\section{Conclusion}
In this paper, we completely determined $\kappa^2(T_kG_n)$, where $T_kG_n$ is a Cayley graph generated by a $k$-tree. As $k$-tree is a generalization of tree and the complete graph is also a $n-1$-tree, the result in this paper
generalizes the result in~\cite{eddie3} by Cheng and Lipt\'ak  on Cayley graphs generated by transposition trees (independently by Yang and Meng~\cite{yang2}), and Wang's result on the complete transposition graphs~\cite{wang}. One can see that the size of the neighborhood of a shortest cycle of a graph usually provides an upper bound for $\kappa^2$ of the graphs. The proof in this paper
shows that $\kappa^2(T_kG_n)$ reaches this upper bound. In this sense,  $T_kG_n$  is optimal with respect to  $\kappa^2$, and thus has high fault tolerance
measured by  $\kappa^2$.

Currently, there are only several results on  $\kappa^2$  for some special Cayley graphs. The problems on determining general transposition networks is still open (Cayley graphs generated by an arbitrary set of transpositions.). Combining Theorem~\ref{uc} and the main result in this paper, we leave a problem below.

\begin{conjecture}
Assume that $G(\mathcal{T})$ is connected. Let $G_n$ the Cayley graph generated by  $G(\mathcal{T})$. Then the following holds.

$(1)$ If  $G(\mathcal{T})$ contains a triangle, then  $\kappa^2(G_n)=4|E(G(\mathcal{T}))|-10$ for $n\geq5$;

$(1)$ If  $G(\mathcal{T})$ contains no triangles, then  $\kappa^2(G_n)=4|E(G(\mathcal{T}))|-8$ for $n\geq4$.
\end{conjecture}


\begin{thebibliography}{13} \small

\bibitem{akers}S.B. Akers, D. Harel, B. Krishnamurthy, The star graph: An attractive alternative to the $n$-cube, Proc. Int. Conf. Parallel Process. (1987) 393-400.

\bibitem{beineke1} L.W. Beineke, R.E. Pippert, The number of labeled $k$-dimensional trees, Journal of Combinatorial Theory 6 (1969) 200-205.

\bibitem{beineke2} L.W. Beineke, R.E. Pippert, Properties and characterizations of $k$-trees, Mathematika 18 (1971) 141-151.

 \bibitem{chang}   N. Chang, S. Hsieh, {2, 3}-Extraconnectivities of hypercube-like networks, Journal of Computer and System Sciences 79 (2013) 669-688.

\bibitem{chou} Z. Chou, C. Hsu, J. Sheu, Bubble-sort star graphs: a new interconnection network, in: International Conference on Parallel and Distributed
Systems, 1996, pp. 41-48.

\bibitem{eddie1} E. Cheng, L. Lipt\'ak, Fault resiliency of Cayley graphs generated by transpositions, International Journal of Foundations of Computer Science 18 (2007)
1005-1022.
\bibitem{eddie2}  E. Cheng, L. Lipt\'ak, Linearly many faults in Cayley graphs generated by transposition trees, Information Sciences 177 (2007) 4877-4882.
\bibitem{eddie3}  E. Cheng, L. Lipt\'ak, A kind of conditional vertex connectivity of Cayley graphs generated by transposition trees, Congressus Numerantium 199 (2009)
167-173.
\bibitem{eddie4}E. Cheng, L. Lipt\'ak, W. Yang, Z. Zhang, X. Guo, A kind of conditional vertex connectivity of Cayley graphs generated
by 2-trees, Information Sciences 181 (2011) 4300-4308.

\bibitem{esfa} A.H. Esfahanian, Generalized Measure of Fault Tolerance with Application to $N$-cube
Networks,\quad IEEE Trans. Comput. 38 (1989) 1586-1591.

\bibitem{fabrega}J. F\`abrega and M.A. Fiol, On the extraconnectivity
of graphs,\quad Discrete Mathematics, 155 (1996) 49-57.

\bibitem{harary}F. Harary, Conditional connectivity,\quad Networks,
13 (1983) 346-357.

\bibitem{heydemann} M. Heydemann, Cayley graphs and interconnection networks, in: G. Hahn, G. Sabidussi (Eds.), Graph Symmetry, Kluwer Academic Publishers.,
Netherlands, 1997, pp. 167-224.



\bibitem{huyang} S.C. Hu, C.B. Yang, Fault tolerance on star graphs, Proceedings of the First Aizu International Symposium on Parallel Algorithms/Architecture Synthesis
(1995) 176-182.

\bibitem{jwo} S. Lakshmivarahan, J. Jwo, S.K. Dhall, Symmetry in interconnection networks based on Cayley graphs of permutation groups: a survey, Parallel
Computing 19 (1993) 361-407.

\bibitem{lai}P. Lai, J. Tan, C. Chang, and L. Hsu, Conditional Diagnosability Measures for Large Multiprocessor
Systems,\quad  IEEE Tans. Comput.  54 (2005) 165-175.

\bibitem{latifi}S. Latifi, M. Hegde and
M. Naraghi-Pour, Conditional connectivity measures for large
multiprocessor systems,\quad IEEE  Trans.  Comput. 43 (1994)
218-222.

\bibitem{latifi2} Latifi, S. and P.K. Srimani, Transposition networks as a class of fault-tolerant robust networks, IEEE Trans. Parallel
Distrib. Sys. 45(2) 1996 230-238.

\bibitem{oh}A.D. Oh and H. Choi, Generalized Measures of Fault Tolerance in $n$-Cube
Networks,\quad IEEE Trans. Parallel
Distrib. Sys.  4 (1993) 702-703.

\bibitem{tindell} R. Tindell, Connectivity of Cayley digraphs, in: Ding-Zhu Du, D. Frank Hsu (Eds.), Combinatorial Network Theory, Kluwer Academic Publishers.,
Netherlands, 1996, pp. 41-64.


\bibitem{wan}M. Wan and Z. Zhang, A kind of conditional vertex connectivity of star graphs. \quad Applied Mathematics Letters 22 (2009) 264-267.

\bibitem{wang} G. Wang, H. Shi, F. Hou, Y. Bai, Some condition vertex connectivities of complete-transposition graphs, Information science 295 (2015) 536-543.


\bibitem{wu}J. Wu and G. Guo, Fault Tolerance Measures for $m$-Ary $n$-Dimensional Hypercubes Based
on Forbidden Faulty Sets,\quad IEEE Trans. Comput. 47 (1998)
888-893.

\bibitem{yang}W. Yang and J. Meng, Extraconnectivity of Hypercubes,\quad Applied Mathematics Letters 22 (2009) 887-891.

\bibitem{yang1} W. Yang and J. Meng, Extraconnectivity of Hypercubes
(II), Australasian Journal of Combinatorics 47 (2010) 189-195.

\bibitem{yang2} W. Yang, H.Z. Li, J.X. Meng, Conditional connectivity of Cayley graphs generated by transposition trees, Information Processing Letters 110 (2010)
1027-1030.

\bibitem{yang3} W. Yang, J.X. Meng, Generalized measures of fault tolerance in hypercube networks, Applied Mathematics Letters 25 (2012) 1335-1339.


\bibitem{zhang1} Z. Zhang, W. Xiong, W-H. Yang, A kind of conditional fault tolerance of alternating group graphs, Information Processing Letters 110 (2010) 998-1002.

\bibitem{zhang2}X. Yu, X. Huang, Z. Zhang,  A kind of conditional connectivity of Cayley graphs generated
by unicyclic graphs, Information Sciences 243 (2013) 86-94.


\end{thebibliography}
\end{document}